\documentclass{amsart}
\usepackage{amsmath}
\usepackage{mathrsfs}
\usepackage{graphicx,url,xy,amssymb,booktabs}
\usepackage[noadjust]{cite}
\usepackage{xspace}

\newcommand{\ZZ}{\mathbb Z}
\newcommand{\QQ}{\mathbb Q}
\newcommand{\CC}{\mathbb C}
\newcommand{\OO}{\mathcal O}

\newcommand{\HH}{\mathbb{H}}
\newcommand{\sd}{{\sqrt{d}}}

\DeclareMathOperator{\GL}{GL}

\DeclareMathOperator{\Norm}{N}
\newcommand{\norm}{\Norm_{F/\QQ}}

\renewcommand{\Re}{\operatorname{Re}}
\newcommand{\Vor}{Vorono\"{\i}\xspace}

\newcommand{\vect}[1]{\begin{bmatrix} #1 \end{bmatrix}}
\newcommand{\mat}[1]{\begin{bmatrix} #1 \end{bmatrix}}%

\theoremstyle{plain}
\newtheorem{thm}{Theorem}[section]

\theoremstyle{definition}
\newtheorem{defn}[thm]{Definition}

\begin{document}


\title{Hyperbolic tessellations associated to Bianchi groups}
\author{Dan Yasaki}
\address{Department of Mathematics and Statistics\\146 Petty Building\\ 
University of North Carolina at Greensboro\\Greensboro, NC 27412}
\email{d\_yasaki@uncg.edu}
\date{}
\thanks{Partially supported by UNCG New Faculty Grant.}
\thanks{The original manuscript was prepared with the \AmS-\LaTeX\ macro
system.}

\begin{abstract}
Let $F/\QQ$ be number field.  The space of positive definite binary Hermitian 
forms over $F$ form an open cone in a real vector space.  There is a natural 
decomposition of this cone into subcones, which descend give rise to 
hyperbolic tessellations of 3-dimensional hyperbolic space by ideal polytopes.
 We compute the structure of these polytopes for a range of imaginary quadratic 
fields.  
\end{abstract}
\maketitle

\bibliographystyle{alpha} 
\begin{section}{Introduction}
Let $F/\QQ$ be a number field.  The space of positive definite binary 
Hermitian forms over $F$ form an open cone in a real vector space.  There is 
a natural decomposition of this cone into polyhedral cones corresponding to 
the facets of the \Vor polyhedron \cite{Gmod,koecher,A}.  This has been 
computationally explored for real quadratic fields in \cite{Ong, ants} and 
the cyclotomic field $\QQ(\zeta_5)$ in \cite{Yascyclotomic}.  

For $F$ an imaginary quadratic field, the polyhedral cones give rise to ideal 
polytopes in $\HH_3$, 3-dimensional hyperbolic space.  In work of Cremona and 
his students \cite{Whthesis,Bythesis,Lithesis,Cr,CrWh}, analogous polytopes 
have already been computed for class number 1 fields imaginary quadratic fields 
as well as a few fields with class number 2 and 3 using 
different methods.  The structure of the polytopes was used to compute Hecke 
operators on modular forms for the Bianchi groups over those fields.  These 
polytopes were used by Goncharov \cite{Gon.euler} in his study of Euler 
complexes on modular curves.  The data of the polytope and stabilizer could
also be used  to give explicit presentations of $\GL_2(\OO)$.  We will examine 
this in a later project.

 In this paper, we investigate the structure of these ideal polytopes for a 
much larger range of 
imaginary quadratic fields.   Our approach and implementation works for general
imaginary quadratic fields, but we restrict the range to ease the computation.
   Specifically, we compute the ideal polytope 
classes for all imaginary quadratic fields of class number 1 and 2, as well 
as some fields of higher class number with small discriminant.  Specifically, 
we compute the ideal polytopes for  the fields $\QQ(\sqrt{d})$ for square-free 
$d$, where 
\[-d \in \{1 - 100,115,123,163,187,235,267,403,427\}.\]
There is no theoretical obstruction to computing these tessellations for 
higher class number and higher discriminant.  

The structure of the paper is as follows.  We set the notation for the 
quadratic fields and Hermitian forms in Section~\ref{sec:notation}.
The implementation is described in Section~\ref{sec:implementation}.  
Finally, in Section~\ref{sec:results}, we summarize some of the data collected 
so far.  

I would like to thank John Cremona for helpful conversations at the beginning 
of this project, and Paul Gunnells for introducing me to these techniques.
I thank Sebastian Pauli for his advice on the computation, and Carlos Nicholas 
for his help with the polytopes.  Finally, I thank Steve Donnelly for helpful 
discussions and the Magma Group at the University of Sydney for their 
hospitality during a visit, in which part of this research was completed.
\end{section}
\begin{section}{Notation and background}\label{sec:notation}
Let $F = \QQ(\sd ) \subset \CC$ be an imaginary quadratic number field.  We 
always take $d<0$ to be a square-free integer.  Let $\OO \subset F$ denote 
the ring of integers in $F$.  Then $\OO$ has a $\ZZ$-basis consisting of $1$ 
and $\omega$, where 
\[\omega = \begin{cases}
\frac{1+\sd}{2} & \text{if $d \equiv 1 \bmod{4}$,}\\
\sd & \text{if $d \equiv 2,3 \bmod{4}$.}
\end{cases}\]
Then \cite{Stark1, Baker2}
$F$ has class number $h_F=1$ if \[-d \in\{1,2,3,7,11,19,43,67,163\}\] and 
$h_F = 2$ if 
\[-d \in \{5,6,10,13,15,22,35,37,51,58,91,115,123,187,235,267,403,427\}.\]

Let $\bar{\cdot}$ denote complex conjugation, the nontrivial Galois 
automorphism of $F$.

\begin{defn}
A \emph{binary Hermitian form over $F$} is a map $\phi:F^2 \to \QQ$ of the form 
\[\phi(x,y) = a x \bar{x} + b x \bar{y} + \bar{b}\bar{x} y + c y \bar{y},\]
where $a,c \in \QQ$ and $b\in F$ such that $\phi$ is positive definite.
\end{defn}
By choosing a $\QQ$-basis for $F$, $\phi$ can be viewed as a quadratic form over $\QQ$.  In particular, it follows that $\phi(\OO^2)$ is discrete in $\QQ$.  

\begin{defn}  
The \emph{minimum of $\phi$} is 
\[m(\phi)=\inf_{v \in \OO^2 \setminus \{ 0\}} \phi(v).\]
A vector $v\in \OO^2$ is \emph{minimal vector} for $\phi$ if 
$\phi(v)=m(\phi)$.  The set of minimal vectors for $\phi$ is denoted 
$M(\phi)$. 
\end{defn}

\begin{defn}
A Hermitian form over $F$ is \emph{perfect} if it is uniquely determined by
 $M(\phi)$ and $m(\phi)$. 
\end{defn}

\end{section}

\begin{section}{Implementation}\label{sec:implementation}
The space of positive definite binary Hermitian forms over $F$ form an 
open cone in a real vector space.   There is 
a natural decomposition of this cone into polyhedral cones corresponding to 
the facets of the \Vor polyhedron $\Pi$ \cite{Gmod,koecher,A}.  The 
top-dimensional cones of this decomposition correspond to perfect forms and 
descend to ideal polytopes in $\HH_3$, 3-dimensional hyperbolic space. 

There is an algorithm \cite{Gmod} to compute the $\GL_2(\OO)$-equivalency 
classes of perfect forms.  The algorithm uses linear 
algebra and convex geometry, but requires an initial input of a perfect form.  
To this end, we describe the method that was used to compute an initial 
perfect form.  

A perfect form $\phi$ is uniquely determined by its minimum $m(\phi)$ and 
set of minimal vectors $M(\phi)$.  By scaling, we can assume $m(\phi) = 1$.
  Since each minimal vector defines a linear equation in $V$, and $V$ is 
4-dimensional, generically 4 minimal vectors will uniquely determine $\phi$.
  Note that this does not imply that $\#M(\phi) = 4$.  Indeed in many 
examples $M(\phi) > 4$.

For each field $F = \QQ(\sqrt{-d})$, we need only to find a single perfect 
form to begin the algorithm.  Thus we limit our search to a particular family 
of quadratic forms.  Specifically, let $S_0 \subset V$ be the subset of 
quadratic forms $\phi$ such that 
\[\left\{ \vect{1\\0}, \vect{0\\1},\vect{1\\1}\right\} \subseteq M(\phi).\]
For $\phi \in S_0$, the Hermitian matrix $A_\phi$ associated to $\phi$ must 
have the form 
\[A_\phi = \mat{1 & \beta\\\bar{\beta} &  1}, \quad 
\text{where $\beta \in F$ with $\Re(\beta)  = -\frac{1}{2}$}.\] 

If $\phi \in S_0$ and $\phi$ has an additional minimal vector 
$\vect{a\\b} \in \OO^2$, then 
\begin{equation}\label{eq:abmin}
\beta = -\frac{1}{2}+
\left(\frac{1-{a_1}^{2}+{a_{{2}}}^{2}d+a_{{1}}b_{{1}}-a_{{2}}db_{{2}
}-{b_{{1}}}^{2}+{b_{{2}}}^{2}d}{2\,da_{{1}}b_{{2}}-2\,da_{{2}}b_{{1}}}
\right) \sqrt{d},
\end{equation}
where $a = a_1+a_2\sqrt{d}$ and $b=b_1 +b_2\sqrt{d}$.  Since $\phi$ is 
positive definite, we must have $\beta \bar{\beta} < 1$.  Combined with 
\eqref{eq:abmin}, this implies
\begin{equation} \label{eq:abbound}
-{\frac { \left( 1-{a_{{1}}}^{2}+{a_{{2}}}^{2}d+a_{{1}}b_{{1}}-a_{{2}}
db_{{2}}-{b_{{1}}}^{2}+{b_{{2}}}^{2}d \right) ^{2}d}{ \left( 2\,da_{{1
}}b_{{2}}-2\,da_{{2}}b_{{1}} \right) ^{2}}}<\frac{3}{4}.
\end{equation}

By reduction theory, the values $\norm(a)$, $\norm(b)$, and $\norm(b-a)$ 
are bounded above by a constant depending upon $d$.  Thus we implement a 
brute force search over $a,b \in \OO$ beginning at $0$ and moving out.  When 
a vector $\vect{a\\b}$ is found satisfying \eqref{eq:abbound}, we check that 
the corresponding form $\phi$ satisfies 
\[\left\{\vect{1\\0}, \vect{0\\1}, \vect{1\\1}, \vect{a\\b}\right\} 
\subseteq M(\phi).\]
This corresponds to a ideal polytope whose vertices contain 
$\{\infty, 0, 1, \frac{a}{b}\}$.

Once the initial form is found, we implement the algorithm of \cite{Gmod} 
to find all the perfect forms and the full structure of the \Vor polyhedron
 in \verb+Magma+ \cite{magma}.
\end{section}
\begin{section}{Results}\label{sec:results}
In this section we collect the results of the computations of the 
$\GL_2(\OO)$-conjugacy classes of the ideal \Vor polytopes.  

\subsection{Example: $d = -14$}
Let $F = \QQ(\sqrt{-14})$.  Then $F$ has class number 4 and ring of integers $\OO = \ZZ[\omega]$, where $\omega = \sqrt{-14}$.  There are 9 $\GL_2(\OO)$-classes of polytopes which are of 3 combinatorial 
types.  There are 3 triangular prisms with cuspidal vertices
\begin{align*}
P_1 &= \left\{\infty, 1, \frac{5 + 2\omega}{9}, 
\frac{2 + \omega}{4}, \frac{4 
+ 2\omega}{9}, 0\right\}\\
P_2&=\left\{\frac{11 + 4\omega}{23}, 1, \frac{5 + 2\omega}{9}, \frac{4 + 
2\omega}{9}, \frac{12 + 4\omega}{23}, 0\right\}, \quad \text{and}\\
P_3 & = \quad\left\{\frac{8 + 5\omega}{23}, \frac{2 + \omega}{5}, 
\frac{1 + \omega}{5}, 
\frac{2 + \omega}{6}, \frac{3 + 2\omega}{10}, \frac{7 + 
4\omega}{21}\right\},\end{align*}
and 5 tetrahedra with cuspidal vertices
\begin{align*}
T_1& =\left\{ \frac{11 + 4\omega}{23}, \frac{2 + \omega}{5}, \frac{4 + 
2\omega}{9}, 0\right\},\\
T_2 & = \left\{ 1, \frac{5 + 2\omega}{9}, \frac{3 + \omega}{5}, \frac{12 + 
4\omega}{23}\right\},\\
T_3 & =\left\{ \frac{11 + 4\omega}{23}, \frac{2 + \omega}{5}, \frac{2 + 
\omega}{6}, 0\right\},\\
T_4 & = \left\{ \frac{8 + 5\omega}{23}, \frac{2 + \omega}{5}, \frac{4 + 
2\omega}{9}, 0\right\},\quad \text{and}\\
T_5 &= \left\{ \frac{4 + \omega}{6}, 1, \frac{3 + \omega}{5}, \frac{12 + 
4\omega}{23}\right\},
\end{align*}
and a square pyramid with cuspidal vertices
\[S=\left\{\frac{8 + 5\omega}{23}, \frac{2 + \omega}{5}, \frac{1 + \omega}{5}, 
\frac{2 + \omega}{6}, 0\right\}.\]
Given the cuspidal vertices, one can easily compute the stabilizers of each polytope.  The stabilizers are all cyclic in this case.  For each stabilizer, we compute a generator.  The results are given in Table~\ref{tab:stabilizers}.
\begin{table}
\caption{Stabilizer groups of \Vor ideal polytopes for $\QQ(\sqrt{-14})$}
\label{tab:stabilizers}
\begin{tabular}{c c c}\toprule
Polytope & Stabilizer & Generator\\\midrule
$P_1$ & $C_6$ & $\mat{1 & -1\\1 & 0}$\\
$P_2$ & $C_2$ & $\mat{-1 & 0\\0 & -1}$\\
$P_3$ & $C_4$ & $\mat{ \omega + 1 &-\omega + 6\\2 &-\omega - 1}$\\
$T_1$ & $C_2$ & $\mat{-1 & 0\\0 & -1}$\\
$T_2$ & $C_2$ & $\mat{-1 & 0\\0 & -1}$\\
$T_3$ & $C_2$ & $\mat{-1 & 0\\0 & -1}$\\
$T_4$ & $C_2$ & $\mat{-1 & 0\\0 & -1}$\\
$T_5$ & $C_2$ & $\mat{-1 & 0\\0 & -1}$\\
$S$ & $C_2$ & $\mat{-1 & 0\\0 & -1}$\\\bottomrule
\end{tabular}
\end{table}
\subsection{Summary}
We compute the \Vor polytopes for all imaginary quadratic number fields $F = \QQ(\sqrt{d})$ with class number $1$ and $2$  as well as higher class number for $d >-100$.   Although there is no reason an arbitrary convex 3-dimensional polytope could not arise, in all of these cases only 8 combinatorial types show up.  We give the names and $F$-vector ([\#vertices, \#edges, \#faces]) for each in Table~\ref{tab:types}.  We also note that the triangular dipyramid shows up in this range much less frequently than the other polytopes.  
\begin{table}
\caption{Combinatorial types of ideal polytopes that occur in this range.}
\label{tab:types}
\begin{tabular}{c c c}
\toprule
polytope & $F$-vector & picture\\
\midrule
 tetrahedron &$[4,6,4]$& \includegraphics[width=0.06\textwidth]{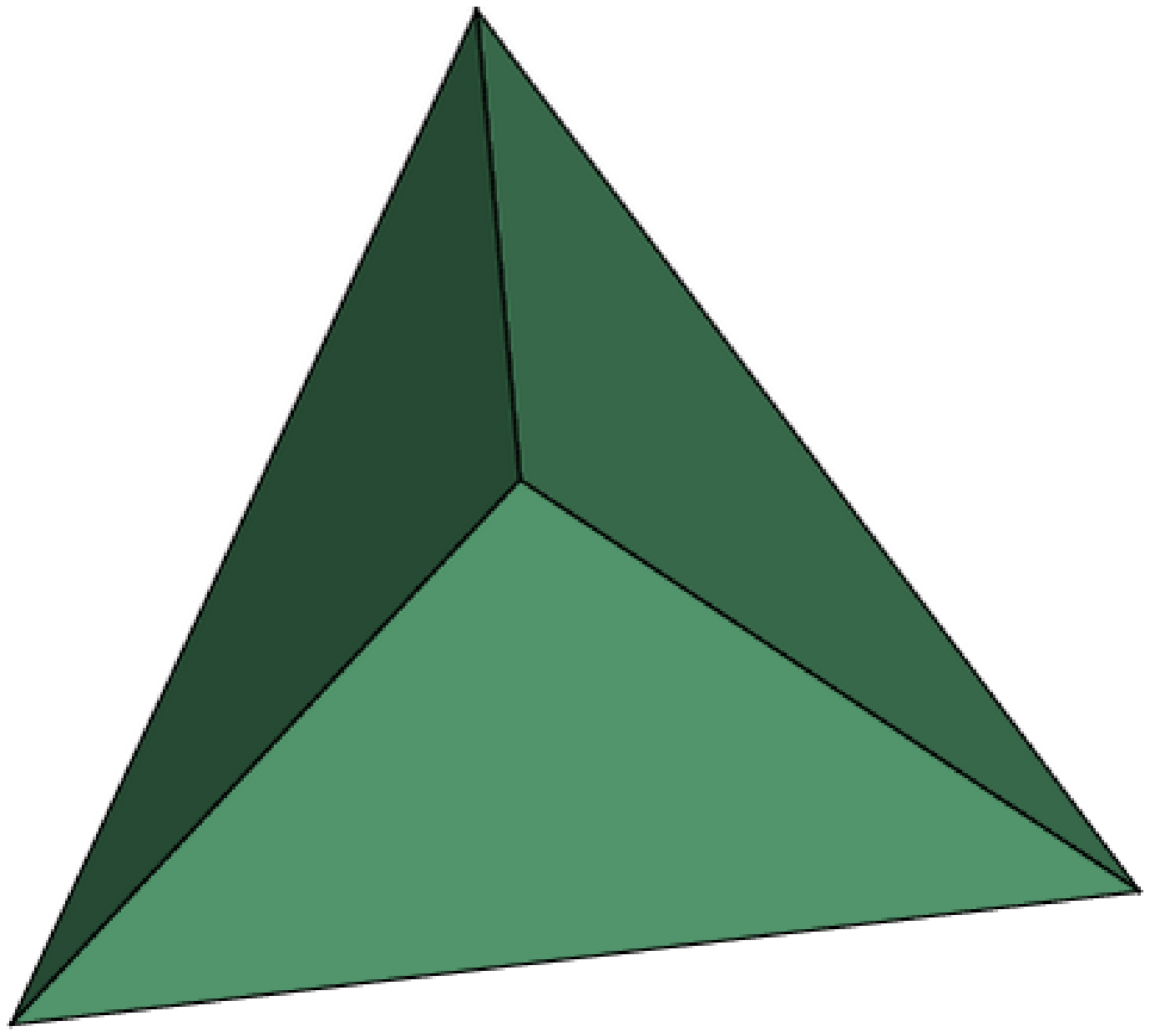}\\
 octahedron& $[6,12,8]$ &\includegraphics[width=0.06\textwidth]{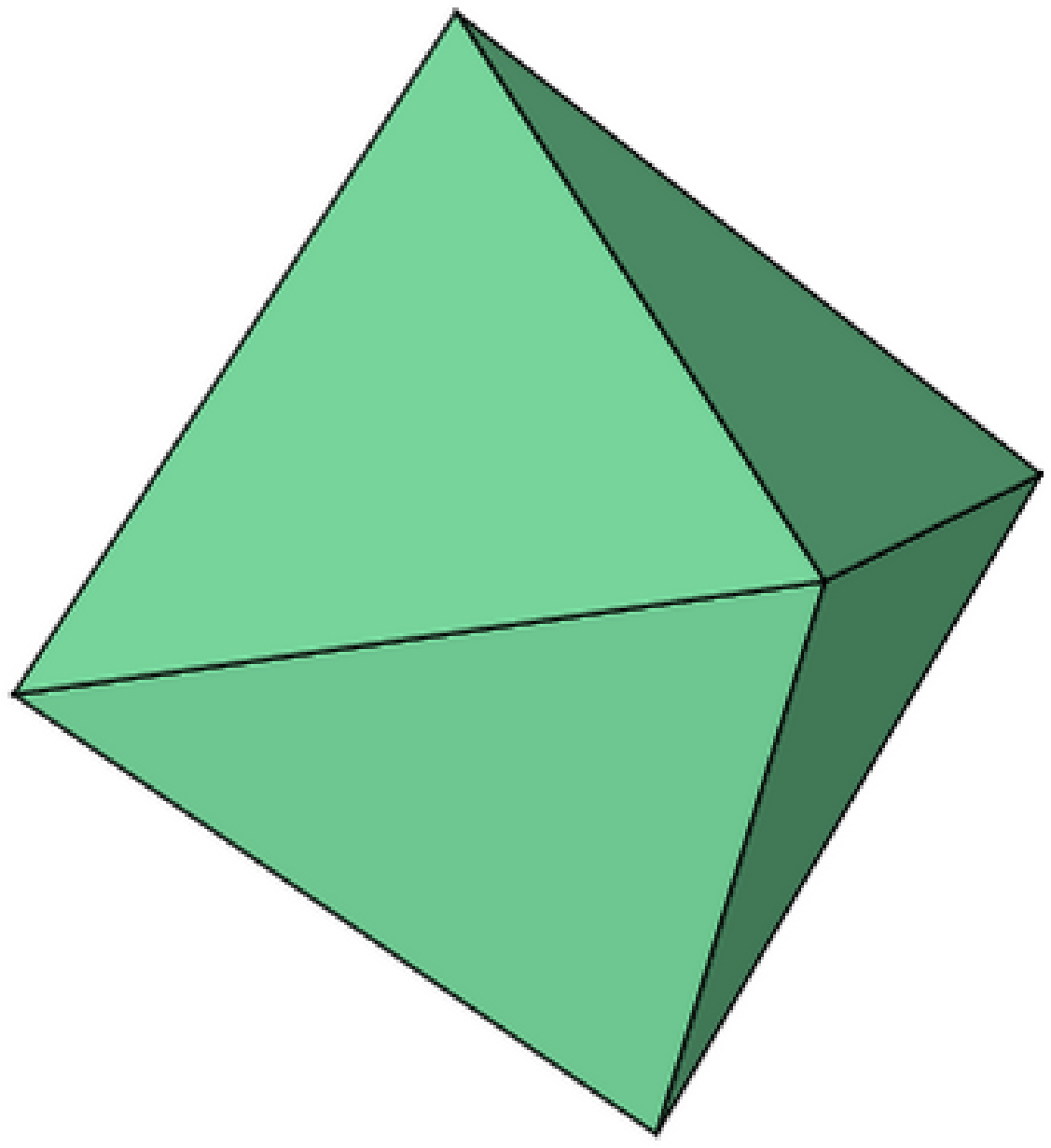}\\
cuboctahedron& $[12,24,14]$&  \includegraphics[width=0.06\textwidth]{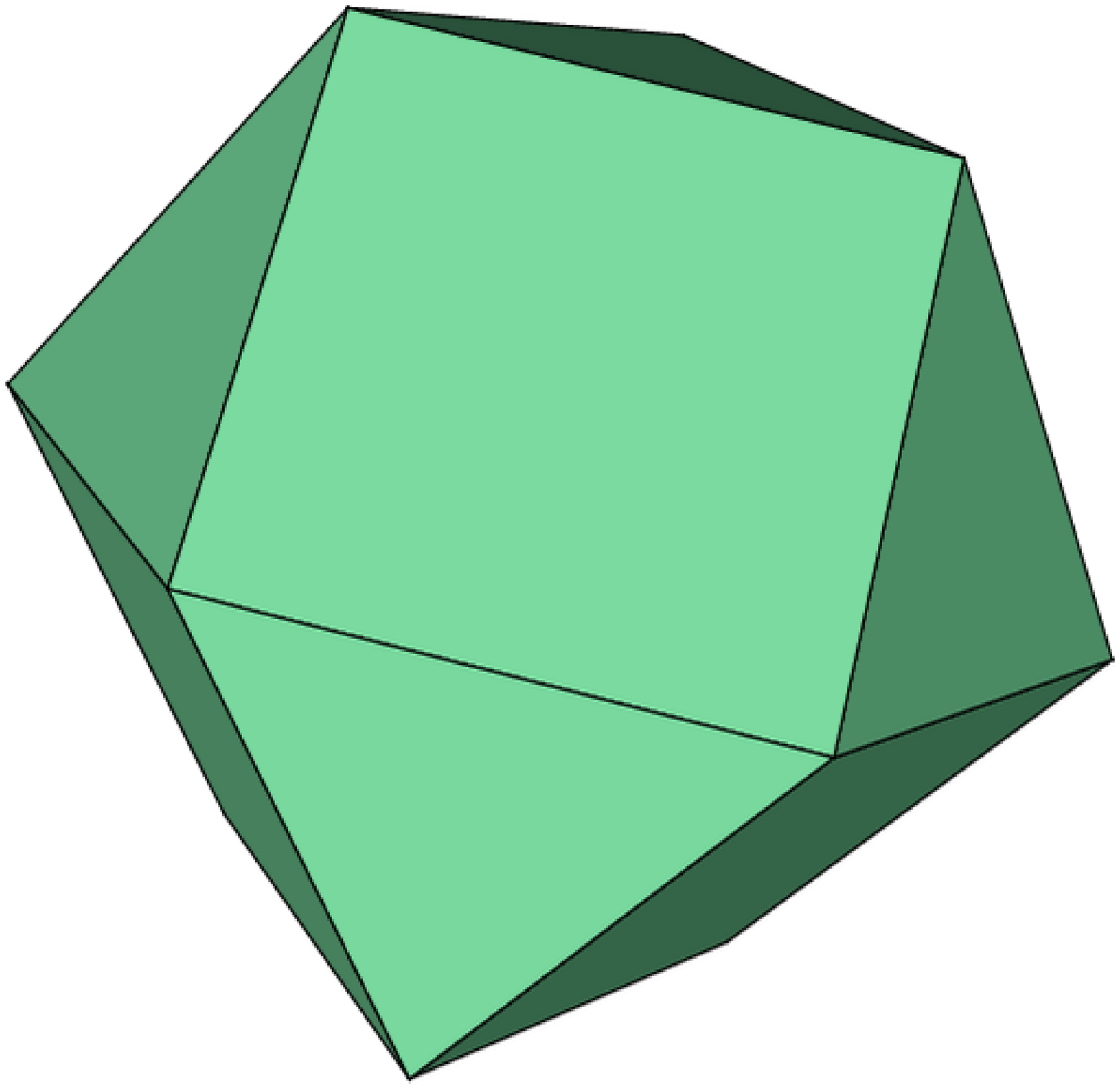}\\
 triangular prism& $[6,9,5]$&\includegraphics[width=0.06\textwidth]{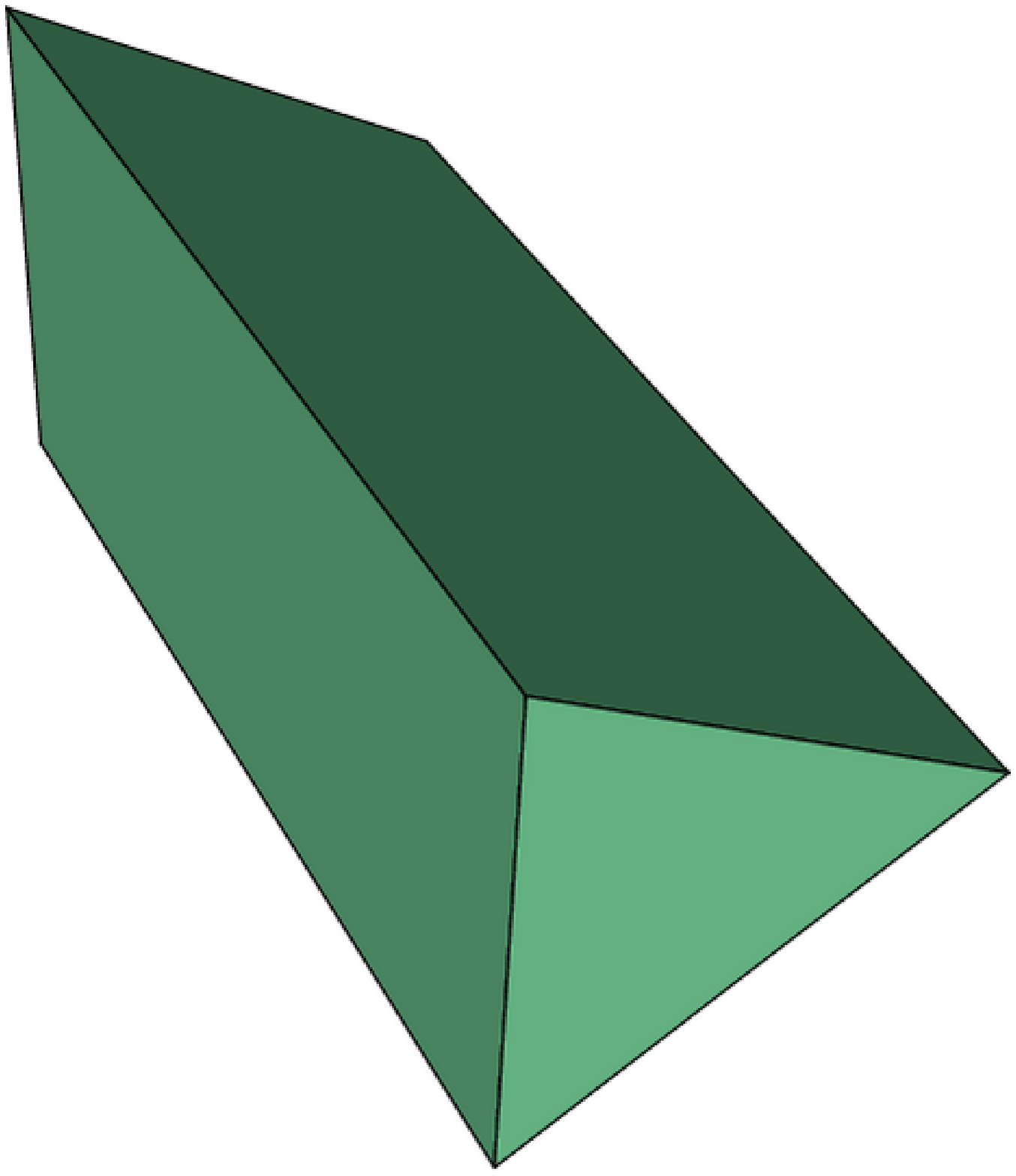}\\
 hexagonal cap& $[9,15,8]$ &\includegraphics[width=0.06\textwidth]{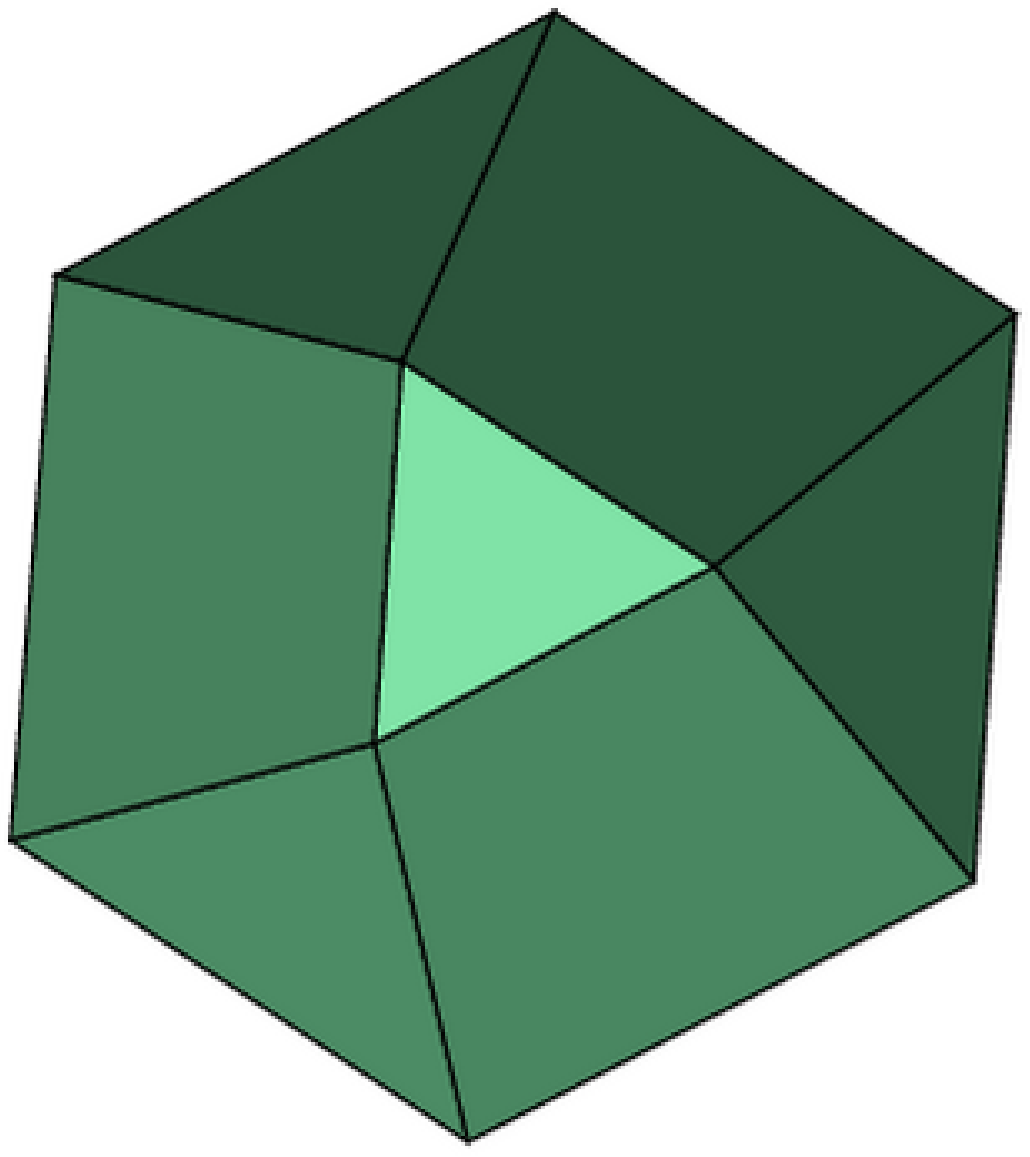}\\
square pyramid& $[5,8,5]$ &\includegraphics[width=0.06\textwidth]{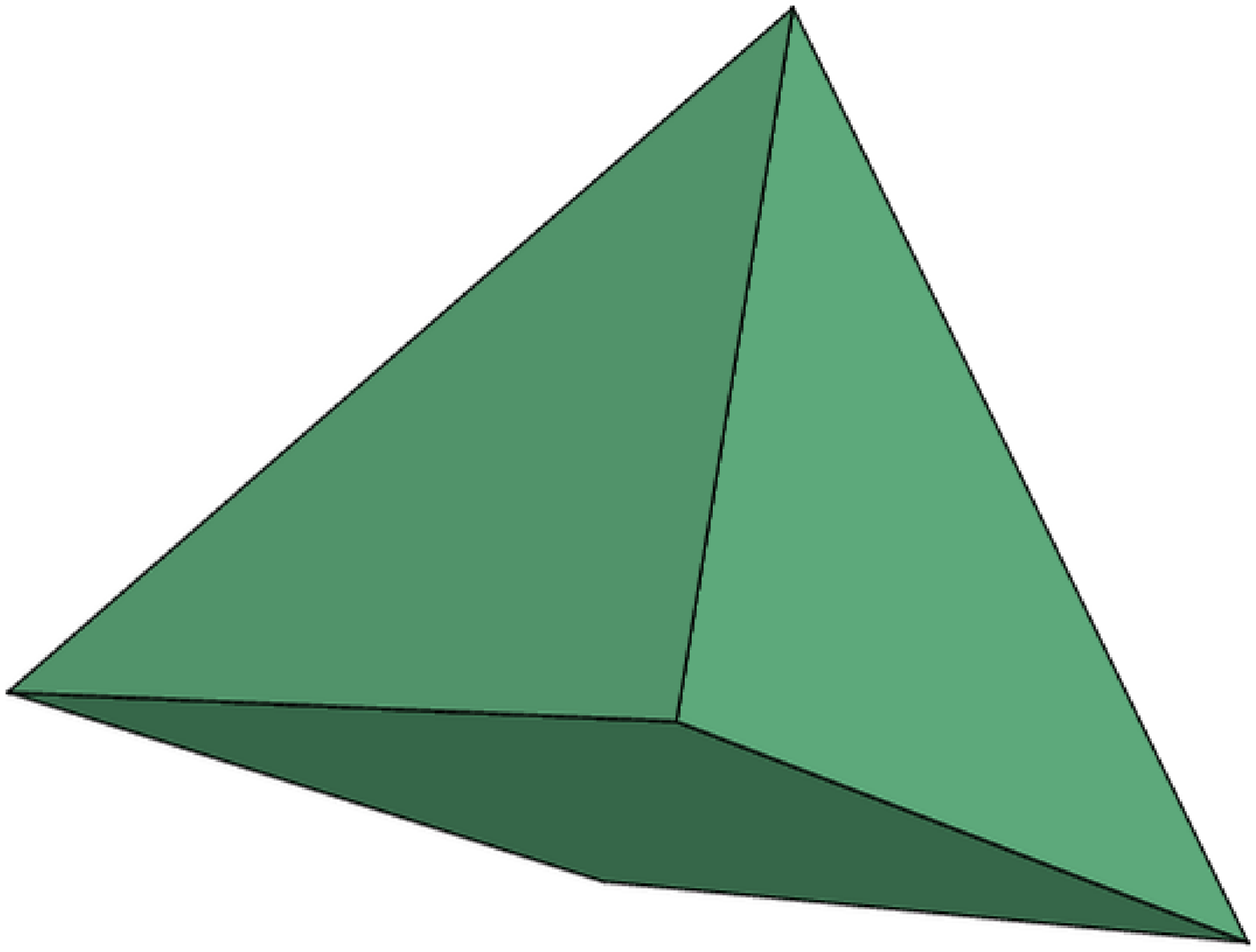}\\
 truncated tetrahedron &$[12,18,8]$& \includegraphics[width=0.06\textwidth]{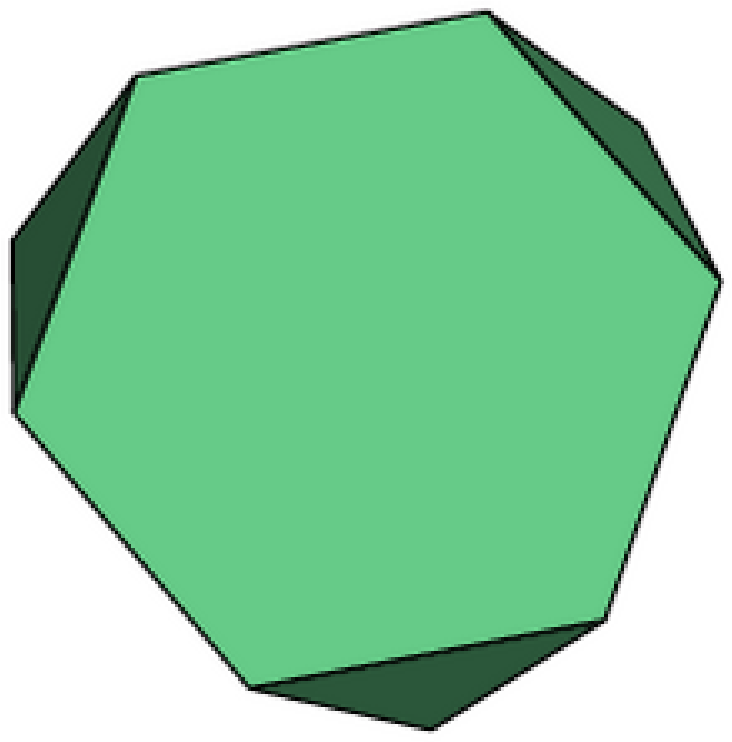}\\
 triangular dipyramid &$[5,9,6]$& 
\includegraphics[width=0.05\textwidth]{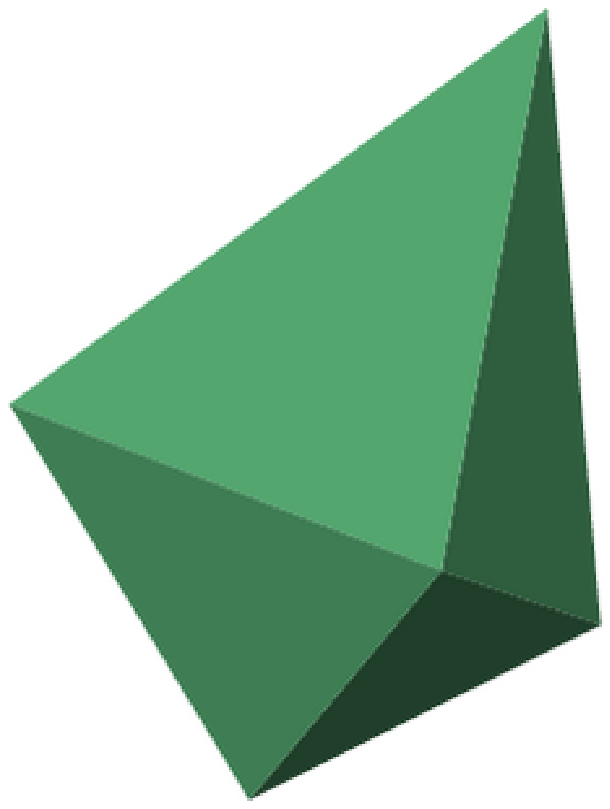}\\
\bottomrule
\end{tabular}
\end{table} 

In Table~\ref{tab:summary12}, we give the number of $\GL_2(\OO)$-classes of each polytope type for $F$ with class number 1 or 2.  In Table~\ref{tab:summaryrest}, we give the number of $\GL_2(\OO)$-classes of each polytope type for the remaining imaginary quadratic fields with $d >-100$.

\begin{table}
\caption{$\GL_2(\OO)$-classes of \Vor ideal polytopes for class number 1 and 2.}
\label{tab:summary12}
\begin{tabular}{cccccccccc}\toprule
  
  $h_F$&$d$&\includegraphics[width=0.06\textwidth]{tetrahedron}&
\includegraphics[width=0.06\textwidth]{octahedron}&
\includegraphics[width=0.06\textwidth]{cuboctahedron}&
\includegraphics[width=0.06\textwidth]{triangularprism}&
\includegraphics[width=0.06\textwidth]{9-15-8}&
\includegraphics[width=0.06\textwidth]{squarepyramid}&
\includegraphics[width=0.06\textwidth]{truncatedtetrahedron}&
\includegraphics[width=0.05\textwidth]{triangulardipyramid}\\
\midrule

1 & $-1$&0&1&0&0&0&0&0&0\\
1 &$-2$&0&0&1&0&0&0&0&0\\
1 &$-3$&1&0&0&0&0&0&0&0\\
1 &$-7$&0&0&0&1&0&0&0&0\\
1 &$-11$&0&0&0&0&0&0&1&0\\
1 &$-19$& 0&0 &1 &1&0&0&0&0 \\
1 &$-43$ &0  &0  &0  & 2 & 1 &0  & 1&0\\
1 &$-67$ & 0 & 1 & 0 & 2 & 1 & 2 & 1 & 0\\
1 & $-163$ & 11 & 0 & 1 & 8 & 2 & 3 & 0 & 0  \\
\midrule
2 &$-5$&0&0&0&2&0&0&0&0\\  
2 &$-6$&0&0&0&0&1&0&1&0\\
2 &$-10$&0 & 1 &0&1&0& 2&0&0\\
2 &$-13$&1&0 &0 & 3 & 1&1&0&0\\
2 &$-15$ & 1 & 1&0&0&0&0&0&0\\
2 &$-22$ & 5 & 0 & 1 & 4 & 0 & 2 & 0 & 0\\
2 & $-35$ & 3 & 4 & 0 & 1 & 0 & 2 & 0 &0\\
2 & $-37$ & 10 & 0 & 0 & 8 & 1 & 8 & 0 & 0\\
2 & $-51$ & 1 & 0 & 1 & 2 & 1 & 0 & 1 & 0\\
2 & $-58$ & 47 & 0 & 0 & 7 & 2 & 6 & 0 & 0\\
2 & $-91$ & 5 & 1 & 0 & 5 & 0 & 3 & 0 & 0\\
2 & $-115$ & 3 & 1 & 0 & 5 & 2 & 4 & 0 & 0\\
2 & $-123$ & 1 & 1 & 1 & 6 & 3 & 3 & 1 & 0\\
2 & $-187$ & 18 & 1 & 1 & 4 & 1 & 9 & 1 & 0\\
2 & $-235$ & 13 & 1 & 0 & 12 & 4 & 11 & 0 & 0\\
2 & $-267$ & 24 & 1 & 1 & 13 & 5 & 10 & 1 & 0\\
2 & $-403$ & 66 & 1 & 0 & 16 & 2 & 20 & 0 & 2\\
2 & $-427$ & 65 & 2 & 0 & 19 & 4 & 24 & 0 & 0\\
\bottomrule
\end{tabular}
\end{table}

\begin{table}
\caption{$\GL_2(\OO)$-classes of \Vor ideal polytopes with $d>-100$.}
\label{tab:summaryrest}
\begin{tabular}{cccccccccc}\toprule
  
  $h_F$&$d$&
\includegraphics[width=0.06\textwidth]{tetrahedron}&
\includegraphics[width=0.06\textwidth]{octahedron}&
\includegraphics[width=0.06\textwidth]{cuboctahedron}&
\includegraphics[width=0.06\textwidth]{triangularprism}&
\includegraphics[width=0.06\textwidth]{9-15-8}&
\includegraphics[width=0.06\textwidth]{squarepyramid}&
\includegraphics[width=0.06\textwidth]{truncatedtetrahedron}&
\includegraphics[width=0.05\textwidth]{triangulardipyramid}\\
\midrule
3 &$-23$ &  0 &  1   &  0 & 1  & 0  & 1 & 0 & 0  \\
3 & $-31$ & 0 & 0 &  0& 3 & 0 & 1 & 0& 0 \\
3 & $-59$ & 0 &1 &1 &3 &0 &2 &0 &0\\
3 & $-83$ & 6 & 0 &0 &2 &2 &1 &1 &0\\
\midrule
4 &$-14$&5& 0& 0& 3 & 0 & 1&0&0\\
4 &$-17$&5&0 &0 &2&1&3& 1 & 0\\
4 & $-21$ & 8 & 2 & 0 & 2 & 1 & 4 & 0&0\\ 
4 & $-30$ & 6 &0  & 0 & 6 & 4 & 4 & 0&0\\
4 & $-33$ & 9 & 0 & 1 & 8 & 1 & 6 & 1&0\\
4 & $-34$ & 20 &0  &0  & 3 & 1 & 6 & 1&0\\
4 & $-39$ & 1 & 0 & 0 & 3 & 1 & 1 &0 &0\\
4 & $-46$ & 32 & 1 &  0 & 5 & 0 & 9 &  0 &0\\
4 & $-55$ & 5 & 1 & 0 & 2 & 0 & 2 & 0 & 0\\
4 & $-57$ & 33 & 1 & 0 & 10 & 3 & 14 & 2 & 0\\
4 & $-73$ & 57 & 1 & 1 &13 &1 &14 &0 &2\\
4 & $-78$ & 69 &1 &0 &11 &4 &18 &0 &0 \\
4 & $-82$ & 92 &0 &0 &8 &3 &11 &1 &0\\
4 & $-85$ & 56 &0 &0 &17 &0 &28 &0 &0\\
4 & $-93$ & 79 &1 &0 &20 &7 &21 &0 &0\\
4 & $-97$ & 95 &0 &1 &19 &3 &19 &0 &0 \\
\midrule
5 & $-47$ & 5 & 0 &  0 & 1 & 1 & 2 &  0 &0\\
5 & $-79$ & 9 &0 &0 &5 &0 &4 &0 &0 \\
\midrule
6 &$-26$ & 18 & 1 & 0 & 2 & 1 & 4 & 0&0\\
6 &$-29$ & 15 & 0 & 0 & 6 & 0 & 6 & 0&0\\
6 & $-38$ & 33 & 1 &  0& 2 & 1 & 6 & 1&0\\
6 & $-53$ & 45 & 0 & 0 & 7 & 2 & 13 & 0 & 0\\
6 & $-61$ & 41 &1 &0 &11 &1 &16 &0 &0\\
6 & $-87$ & 6 &0 &0 &6 &2 &3 &0 &0 \\
\midrule
7 & $-71$ & 7 &1 &0 &4 &0 &4 &0 &0\\
\midrule
8 & $-41$ & 31 & 0 &  1 & 9 & 0 & 8 &  0 &0\\
8 & $-62$ & 81 &0 &0 &7 &2 &7 &0 &0\\
8 & $-65$ & 69 & 2&0 &9 &0 &19 &0 &0\\
8 & $-66$ & 67 &1 &1 &9 &4 &12 &1 &0\\
8 & $-69$ & 51& 2& 0 &15 &2 &21 &0 &0\\
8 & $-77$ & 81 &1 &0 &9 &2 &26 &0 &0\\
8 & $-94$ & 125 &1 &0 &10 &2 &17 &0 &0\\
8 & $-95$ & 12 &0 &0 &4 &0 &9 &0 &0 \\
\midrule
10 & $-74$ & 105 &1 &0 &9 &1 &12 &0 &0\\
10 & $-86$ & 130 &0 &0 &9 &1 &18 &1 &0\\
\midrule
12 & $-89$ & 136 &0 &0 &14 &1 &21 &1 &0\\
\bottomrule
\end{tabular}

\end{table}
\end{section}

\bibliography{references}    

\begin{thebibliography}{BCP97}

\bibitem[Ash77]{A}
Avner Ash.
\newblock Deformation retracts with lowest possible dimension of arithmetic
  quotients of self-adjoint homogeneous cones.
\newblock {\em Math. Ann.}, 225(1):69--76, 1977.

\bibitem[Bak71]{Baker2}
A.~Baker.
\newblock Imaginary quadratic fields with class number {$2$}.
\newblock {\em Ann. of Math. (2)}, 94:139--152, 1971.

\bibitem[BCP97]{magma}
Wieb Bosma, John Cannon, and Catherine Playoust.
\newblock The {M}agma algebra system. {I}. {T}he user language.
\newblock {\em J. Symbolic Comput.}, 24(3-4):235--265, 1997.
\newblock Computational algebra and number theory (London, 1993).

\bibitem[Byg98]{Bythesis}
J.~Bygott.
\newblock {\em Modular forms and modular symbols over imaginary quadratic
  fields}.
\newblock PhD thesis, Exeter University, 1998.

\bibitem[Cre94]{Cr}
John~E. Cremona.
\newblock Periods of cusp forms and elliptic curves over imaginary quadratic
  fields.
\newblock In {\em Elliptic curves and related topics}, volume~4 of {\em CRM
  Proc. Lecture Notes}, pages 29--44. Amer. Math. Soc., Providence, RI, 1994.

\bibitem[CW94]{CrWh}
J.~E. Cremona and E.~Whitley.
\newblock Periods of cusp forms and elliptic curves over imaginary quadratic
  fields.
\newblock {\em Math. Comp.}, 62(205):407--429, 1994.

\bibitem[Gon08]{Gon.euler}
Alexander~B. Goncharov.
\newblock Euler complexes and geometry of modular varieties.
\newblock {\em Geom. Funct. Anal.}, 17(6):1872--1914, 2008.

\bibitem[Gun99]{Gmod}
Paul~E. Gunnells.
\newblock Modular symbols for {${\bf Q}$}-rank one groups and {V}orono\u\i\
  reduction.
\newblock {\em J. Number Theory}, 75(2):198--219, 1999.

\bibitem[GY08]{ants}
Paul~E. Gunnells and Dan Yasaki.
\newblock Hecke operators and {H}ilbert modular forms.
\newblock In {\em Algorithmic number theory}, volume 5011 of {\em Lecture Notes
  in Comput. Sci.}, pages 387--401. Springer, Berlin, 2008.

\bibitem[Koe60]{koecher}
Max Koecher.
\newblock Beitr\"age zu einer {R}eduktionstheorie in {P}ositivit\"atsbereichen.
  {I}.
\newblock {\em Math. Ann.}, 141:384--432, 1960.

\bibitem[Lin05]{Lithesis}
Mark Lingham.
\newblock {\em Modular forms and elliptic curves over imaginary quadratic
  fields}.
\newblock PhD thesis, University of Nottingham, 2005.

\bibitem[Ong86]{Ong}
Heidrun~E. Ong.
\newblock Perfect quadratic forms over real-quadratic number fields.
\newblock {\em Geom. Dedicata}, 20(1):51--77, 1986.

\bibitem[Sta67]{Stark1}
H.~M. Stark.
\newblock A complete determination of the complex quadratic fields of
  class-number one.
\newblock {\em Michigan Math. J.}, 14:1--27, 1967.

\bibitem[Whi90]{Whthesis}
E.~Whitley.
\newblock {\em Modular symbols and elliptic curves over imaginary quadratic
  number fields}.
\newblock PhD thesis, Exeter University, 1990.

\bibitem[Yas]{Yascyclotomic}
Dan Yasaki.
\newblock Binary hermitian forms over a cyclotomic field.
\newblock to appear in Journal of Algebra.

\end{thebibliography}
\end{document}